\documentclass[12pt,leqno]{amsart}
\usepackage{float}
\usepackage{amssymb,amsmath,amsthm,amscd}
\usepackage{mathrsfs}
\usepackage{enumerate}

\usepackage[all]{xy} 
\CompileMatrices

\numberwithin{equation}{section}

\newcommand{\C}{{\mathbb C}}
\newcommand{\Q}{\mathbb {Q}}

\newcommand{\Z}{{\mathbb Z}}
\newcommand{\B}{{\mathcal{B}}}

\newcommand{\gl}{{\mathfrak{gl}}}

\newcommand{\seteq}{\mathbin{:=}}

\theoremstyle{plain}
\newtheorem{lemma}{Lemma}[section]
\newtheorem{prop}[lemma]{Proposition}
\newtheorem{theorem}[lemma]{Theorem}
\newcommand{\Prop}{\begin{prop}}
\newcommand{\enprop}{\end{prop}}
\newcommand{\Lemma}{\begin{lemma}}
\newcommand{\enlemma}{\end{lemma}}
\newcommand{\Th}{\begin{theorem}}
\newcommand{\enth}{\end{theorem}}
\newtheorem{corollary}[lemma]{Corollary}
\newcommand{\Cor}{\begin{corollary}}
\newcommand{\encor}{\end{corollary}}
\newtheorem{definition}[lemma]{Definition}

\newcommand{\Def}{\begin{definition}}
\newcommand{\edf}{\end{definition}}

\theoremstyle{definition}

\newtheorem{conjecture}[lemma]{Conjecture}

\newcommand{\g}{{\mathfrak{g}}}

\newcommand{\Hom}{\operatorname{Hom}}

\newcommand{\isoto}[1][]%
{{\mathop{\buildrel{\sim}\over\longrightarrow}\limits_{#1}}}
\renewcommand{\hom}{\operatorname{\it \mathscr{H}\kern-.25em om}}

\newcommand{\eq}{\begin{eqnarray}}
\newcommand{\eneq}{\end{eqnarray}}

\newcommand{\eqn}{\begin{eqnarray*}}
\newcommand{\eneqn}{\end{eqnarray*}}

\newenvironment{tenumerate}{
  \begin{enumerate}
  
  }{\end{enumerate}}

\newcommand{\bnum}{\begin{tenumerate}}
\newcommand{\enum}{\end{tenumerate}}
\newenvironment{anumerate}{
  \begin{enumerate}
  
  }{\end{enumerate}}

\newcommand{\bia}{\begin{anumerate}}
\newcommand{\eia}{\end{anumerate}}
\newcommand{\on}{\operatorname}

\newcommand{\bni}{\begin{tenumerate}}
\newcommand{\eni}{\end{tenumerate}}

\newcommand{\QED}{\end{proof}}

\newcommand{\soplus}{\mathop{\mbox{\normalsize$\bigoplus$}}\limits}

\newcommand{\cl}{\colon}

\newcommand{\id}{\on{id}}
\newcommand{\ba}{\begin{array}}
\newcommand{\ea}{\end{array}}

\newcommand{\bi}{\bni}
\newcommand{\ei}{\eni}

\newcommand{\set}[2]{\left\{#1 \mathbin{;} #2 \right\}}

\newcommand{\hs}{\hspace}

\newcommand{\eqsub}{\begin{subequations}\begin{eqnarray}}
\newcommand{\eneqsub}{\end{eqnarray}\end{subequations}}

\newcommand{\ol}{\overline}

\newcommand{\A}{\mathbf{A}}

\renewcommand{\le}{\leqslant}
\renewcommand{\ge}{\geqslant}
\newcommand{\Ind}{{{\on{Ind}}}}

\newcommand{\nc}{\newcommand}
\nc{\la}{\lambda}
\nc{\lam}{\lambda}
\nc{\U}[1][\g]{U_q(#1)}
\nc{\te}{\tilde{e}}
\nc{\tei}{\tilde{e}_i}
\nc{\tf}{\tilde{f}}
\nc{\tfi}{\tilde{f}_i}
\nc{\tU}{\widetilde U_q(\g)}
\nc{\tE}{\tilde{E}}
\nc{\tF}{\tilde{F}}

\nc{\BZ}{{\mathbb{Z}}}
\nc{\al}{\alpha}
\nc{\qs}{{q}}
\nc{\lan}{\langle}
\nc{\ran}{\rangle}
\nc{\re}{{\mathrm{re}}}
\nc{\wt}{\operatorname{wt}}
\nc{\Uf}[1][\g]{U^-_q(#1)}
\nc{\Ue}{U^+_q(\g)}
\nc{\eps}{\varepsilon}
\nc{\vphi}{\varphi}
\nc{\sphi}{\varphi^*}
\nc{\seps}{\varepsilon^*}

\nc{\nn}{\nonumber}

\nc{\vp}{\varpi}

\nc{\cls}{{\operatorname{cl}}}

\nc{\Wt}{{\operatorname{Wt}}}
\nc{\Us}{U'_q(\g)}
\nc{\La}{\Lambda}
\nc{\ro}{{\rm(}}
\nc{\rf}{{\rm)}}
\nc{\norm}{{\mathrm{norm}}}
\nc{\qbox}{\quad\mbox}
\nc{\braid}{{\mathfrak{B}}}
\nc{\Ad}{\operatorname{Ad}}
\nc{\Aut}{\operatorname{Aut}}
\nc{\dt}[1]{\tilde{\tilde #1}}
\nc{\Sn}{S^{{\mathrm{norm}}}}
\nc{\aff}{{\mathrm{aff}}}
\nc{\rk}{{\mathrm{rk}}}
\nc{\tQ}{\widetilde{Q}}
\nc{\tP}{\widetilde{P}}
\nc{\tW}{\widetilde{W}}
\nc{\Dyn}{\mathrm{Dyn}}
\nc{\tD}{\widetilde{\Delta}}
\nc{\height}{{\operatorname{ht}}}
\nc{\bl}{\bigl}
\nc{\br}{\bigr}
\nc{\Hecke}{\mathrm{H}}
\nc{\HA}{\Hecke^{\mathrm{A}}}
\nc{\HB}{\Hecke^{\mathrm{B}}}
\nc{\K}{\mathrm{K}}
\newcommand{\scbul}{{\,\raise1pt\hbox{$\scriptscriptstyle\bullet$}\,}}
\nc{\vac}{{\phi}}
\nc{\Bt}{\B_\theta(\g)}
\nc{\be}{\begin{enumerate}}
\nc{\ee}{\end{enumerate}}
\nc{\low}{{\mathrm{low}}}
\nc{\upper}{{\mathrm{up}}}
\nc{\Zodd}{\Z_{\mathrm{odd}}}
\nc{\Ft}[1][n]{\mathbb{P}\mathrm{ol}_{#1}}
\nc{\Ftf}[1][n]{\widetilde{\mathbb{P}\mathrm{ol}}_{#1}}
\nc{\KA}{\on{K}^{\mathrm{A}}}
\nc{\KB}{\on{K}^{\mathrm{B}}}
\nc{\Res}{\on{Res}}
\nc{\Fc}[1][{n,m}]{\mathbf{F}_{#1}}
\nc{\tphi}{\tilde{\varphi}}
\nc{\CO}{\mathscr{O}}

\begin{document}

\title
{Symmetric crystals and
affine Hecke algebras of type B}
\author{Naoya ENOMOTO}
\author{Masaki KASHIWARA}
\address{Research Institute for Mathematical Sciences,
Kyoto University, Kyoto 606, Japan
}
\thanks{The second author is partially supported by 
Grant-in-Aid for Scientific Research (B) 18340007,
Japan Society for the Promotion of Science.}


\begin{abstract}
The Lascoux-Leclerc-Thibon conjecture, reformulated and solved by
S.{} Ariki, asserts that the K-group of the representations of
the affine Hecke algebras of type A is isomorphic to
the algebra of functions
on the maximal unipotent subgroup of the group associated with a Lie algebra
$\g$ where $\g$ is $\gl_\infty$ or 
the affine Lie algebra $A^{(1)}_\ell$,
and the irreducible representations correspond to the
upper global bases.
In this note, we formulate analogous conjectures for
certain classes of irreducible representations
of affine Hecke algebras of type B.
\end{abstract}

\maketitle

\section{Introduction}
The purpose of this note is to
formulate and explain 
conjectures on certain classes of irreducible representations
of affine Hecke algebras of type B
analogous to the Lascoux-Leclerc-Thibon conjecture (\cite{LLT}),
reformulated and solved by S.{} Ariki,
on affine Hecke algebras of type A.

Let us begin by recalling the Lascoux-Leclerc-Thibon conjecture,
reformulated and solved by S.{} Ariki (\cite{Ariki}).
Let $\HA_n$ be the affine Hecke algebra of type A of degree $n$.
Let $\KA_n$ be the Grothendieck group of
the abelian category of finite-dimensional $\HA_n$-modules,
and $\KA=\oplus_{n\ge0}\KA_n$.
Then it has a structure of Hopf algebra by the restriction and
the induction (cf.\ \S \ref{subsec:affA}).
The set $I=\C^*$ may be regarded as a Dynkin diagram
with $I$ as the set of vertices and with edges between $a\in I$ and $ap_1^2$
(see \eqref{eq:Dy}).
Here $p_1$ is the parameter of 
the affine Hecke algebra usually denoted by $q$.
Let $\g_I$ be the associated Lie algebra, and
$\g_I^-$ the unipotent Lie subalgebra. Let $U_I$ be the group associated
to $\g_I^-$.
Hence $\g_I$ is isomorphic to a direct sum of copies of $A^{(1)}_\ell$
if $p_1^2$ is a primitive $\ell$-th root of unity
and to a direct sum of copies of $\gl_\infty$ if $p_1$ has an infinite order.
Then $\C\otimes \KA$ is isomorphic to the algebra $\CO(U_I)$
of regular functions on $U_I$.
Let $\U[\g_I]$ be the associated quantized enveloping algebra.
Then $\Uf[\g_I]$ has an upper global basis $\{G^{\upper}(b)\}_{b\in B(\infty)}$.
By specializing $\soplus \C[q,q^{-1}]G^\upper(b)$
at $q=1$, we obtain $\CO(U_I)$.
Then the LLT-conjecture says that the
elements associated to irreducible $\HA$-modules
corresponds to the image of the upper global basis.

In this note,
we shall formulate analogous conjectures for affine Hecke algebras of type B.
In the type B case, we have to replace 
$\Uf[\g_I]$ and its upper global basis with
a new object, the symmetric crystals (see \S\;\ref{sec:crystal}).
It is roughly stated as follows.
Let $\HB_n$ be the affine Hecke algebra of type B of degree $n$.
Let $\KB_n$ be the Grothendieck group of
the abelian category of finite-dimensional modules over $\HB_n$,
and $\KB=\oplus_{n\ge0}\KB_n$.
Then $\KB$ has a structure of a Hopf bimodule over $\KA$.
The group $U_I$ has the anti-involution $\theta$ induced by
the involution $a\mapsto a^{-1}$ of $I=\C^*$. Let $U_I^\theta$
be the $\theta$-fixed point set of $U_I$.
Then $\CO(U_I^\theta)$ is a quotient ring of $\CO(U_I)$.
The action of $\CO(U_I)\simeq\C\otimes\KA$ on $\C\otimes\KB$,
in fact, descends to the action of $\CO(U_I^\theta)$.

We introduce 
$V_\theta(\la)$ (see \S\;\ref{sec:crystal}),
a kind of the $q$-analogue of $\CO(U_I^\theta)$.
Our conjecture is then:
\bi
\item
$V_\theta(\la)$ has a crystal basis and a global basis.
\item
$\KB$ is isomorphic to a specialization of $V_\theta(\la)$ at $q=1$
as an $\CO(U_I)$-module, and
the irreducible representations correspond
to the upper global basis of $V_\theta(\la)$ at $q=1$.
\ei
We exclude the representations of $\HB_n$
such that $X_i$ have an eigenvalue $\pm1$ (see \S\;\ref{sec:aff}).

\section{Symmetric crystals}\label{sec:crystal}

In this section, we shall introduce crystals associated with quantum groups
with an involution.

\subsection{Quantized universal enveloping algebras}
We shall recall the quantized universal enveloping algebra
$U_q(\g)$.
Let $I$ be an index set (for simple roots),
and $Q$ the free $\Z$-module with a basis $\{\al_i\}_{i\in I}$.
Let $(\scbul,\scbul)\cl Q\times Q\to\Z$ be
a symmetric bilinear form such that
$(\al_i,\al_i)/2\in\Z_{>0}$ for any $i$ and
$(\al_i^\vee,\al_j)\in\Z_{\le0}$ for $i\not=j$ where
$\al_i^\vee\seteq2\al_i/(\al_i,\al_i)$.
Let $q$ be an indeterminate and set 
$K\seteq\Q(\qs)$.
We define its subrings $\A_0$, $\A_\infty$ and $\A$ as follows.
\eq
&&\ba{rcl}
\A_0&=&\set{f/g}{f(\qs),g(\qs)\in \Q[\qs],\, g(0)\not=0},\\[3pt]
\A_\infty&=&\set{f/g}{f(\qs^{-1}),g(\qs^{-1})\in \Q[\qs^{-1}],\, g(0)\not=0},
\\[3pt]
\A&=&\Q[\qs,\qs^{-1}].
\ea
\eneq

\begin{definition}\label{U_q(g)}
The quantized universal enveloping algebra $U_q(\g)$ is the $K$-algebra
generated by the elements $e_i,f_i$ and invertible
elements $t_i\ (i\in I)$
with the following defining relations.
\begin{enumerate}[{\rm(1)}]

\item The $t_i$'s commute with each other.

\item
$t_je_i\,t_j^{-1}=q^{(\al_j,\al_i)}\,e_i\ $
and $\ t_jf_it_j^{-1}=q^{-(\al_j,\al_i)}f_i\ $
for any $i,j\in I$.

\item\label{even} $\lbrack e_i,f_j\rbrack
=\delta_{ij}\dfrac{t_i-t_i^{-1}}{q_i-q_i^{-1}}$
for $i$, $j\in I$. Here $q_i\seteq q^{(\al_i,\al_i)/2}$.

\item {\rm(}{\em Serre relation}{\rm)} For $i\not= j$,
\begin{eqnarray*}
&&\sum^b_{k=0}(-1)^ke^{(k)}_ie_je^{(b-k)}_i=\sum^b_{k=0}(-1)^kf^{(k)}_i
f_jf_i^{(b-k)}=0.
\end{eqnarray*}
Here $b=1-(\al_i^\vee,\al_j)$ and
\begin{eqnarray*}
e^{(k)}_i=e^k_i/\lbrack k\rbrack_i!\ ,&& f^{(k)}_i=f^k_i/\lbrack k\rbrack_i!\ ,\\
\lbrack k\rbrack_i=(q^k_i-q^{-k}_i)/(q_i-q^{-1}_i)\ ,
&&\lbrack k\rbrack_i!=\lbrack 1\rbrack_i\cdots \lbrack k\rbrack_i\,.
\end{eqnarray*}
\end{enumerate}
\end{definition}

Let us denote by $\Uf$ (resp.{} $\Ue$)
the subalgebra of $\U$ generated by the $f_i$'s (resp.{} the $e_i$'s).
Let us recall the crystal theory of $\Uf$ (\cite{K1}).
Let $e'_i$ and $e^*_i$ be the operators on $\Uf$ defined by
$$[e_i,a]=\dfrac{(e^*_ia)t_i-t_i^{-1}e'_ia}{q_i-q_i^{-1}}\quad(a\in\Uf).$$
Then these operators satisfy the following formula similar to derivations:
$$e_i'(ab)=e_i'(a)b+(\Ad(t_i)a)e_i'b,\quad
e_i^*(ab)=ae_i^*b+(e_i^*a)(\Ad(t_i)b).$$
Then $\Uf$ has a unique symmetric bilinear form $(\scbul,\scbul)$
such that $(1,1)=1$ and
$$(e'_ia,b)=(a,f_ib)\quad\text{for any $a,b\in\Uf$.}$$
It is non-degenerate and satisfies $(e^*_ia,b)=(a,bf_i)$.
The left multiplication of $f_j$ and $e'_i$ have the commutation relation
$$e'_if_j=q^{-(\al_i,\al_j)}f_je'_i+\delta_{ij},$$
and both the $e_i'$'s and the $f_i$'s satisfy the Serre relations.
Since $e'_i$ and $f_i$ satisfy the $q$-boson relation,
any element $a\in\Uf$ can be written uniquely as
$$a=\sum_{n\ge0}f_i^{(n)}a_n\quad\text{with $e'_ia_n=0$.}$$
We define the modified root operators by
\eq
\te_ia=\sum_{n\ge1}f_i^{(n-1)}a_n\quad\text{and}\quad
\tf_ia=\sum_{n\ge0}f_i^{(n+1)}a_n.
\label{eq:modified}
\eneq
Let $L(\infty)$ be the $\A_0$-submodule of $\Uf$ generated by
the $\tf_{i_1}\cdots\tf_{i_\ell}1$'s ($\ell\ge0$, $i_1,\ldots,i_\ell\in I$\,).
Let us set $B(\infty)=
\set{\tf_{i_1}\cdots\tf_{i_\ell}1\bmod \qs L(\infty)}%
{\text{$\ell\ge0$, $i_1,\ldots, i_\ell\in I$}}
\subset L(\infty)/\qs L(\infty)$.
Then we have
\Th
\bi
\item $\tfi L(\infty)\subset L(\infty)$
and $\tei L(\infty)\subset L(\infty)$,
\item
$B(\infty)$ is a basis of $L(\infty)/\qs L(\infty)$,
\item
$\tf_iB(\infty)\subset B(\infty)$ and $\te_iB(\infty)
\subset B(\infty)\sqcup\{0\}$.
\ei
\enth

\subsection{Global bases}\label{subsec:global}
Let $-$ be the automorphism of $K$ sending $\qs$ to $\qs^{-1}$.
Then $\ol{\A_0}$ coincides with $\A_\infty$. 
Let $V$ be a vector space over $K$,
$L_0$ an $A$-submodule of $V$,
$L_\infty$ an $\A_\infty$- submodule, and
$V_\A$ an $\A$-submodule.
Set $E:=L_0\cap L_\infty\cap V_\A$.

\Def[\cite{K1}]
We say that $(L_0,L_\infty,V_\A)$ is {\em balanced}
if each of $L_0$, $L_\infty$ and $V_\A$
generates $V$ as a $K$-vector space,
and if one of the following equivalent conditions is satisfied.
\bnum
\item
$E \to L_0/\qs L_0$ is an isomorphism,
\item
$E \to L_\infty/\qs^{-1}L_\infty$ is an isomorphism,
\item
$(L_0\cap V_\A)\oplus
(\qs^{-1} L_\infty \cap V_\A) \to V_\A$
      is an isomorphism.
\item
$\A_0\otimes_\Q E \to L_0$, $\A_\infty\otimes_\Q E \to L_\infty$,
        $\A\otimes_\Q E \to V_\A$ and $K \otimes_\Q E \to V$
are isomorphisms.
\enum
\edf

Let $-$ be the ring automorphism of $\U$ sending
$\qs$, $t_i$, $e_i$, $f_i$ to $\qs^{-1}$, $t_i^{-1}$, $e_i$, $f_i$.

Let $\U_\A$ be the $\A$-subalgebra of
$\U$ generated by $e_i^{(n)}$, $f_i^{(n)}$
and $t_i$.
Similarly we define
$\Uf_\A$.

\Th
$(L(\infty),L(\infty)^-,\Uf_\A)$ is balanced.
\enth
Let $$G\colon L(\infty)/\qs L(\infty)\isoto E\seteq L(\infty)\cap L(\infty)^-
\cap \Uf_\A$$ 
be the inverse of $E\isoto L(\infty)/\qs L(\infty)$.
Then $\set{G(b)}{b\in B(\infty)}$ forms a basis of $\Uf$.
We call it a (lower) {\em global basis}.
It is first introduced by G.\ Lusztig (\cite{L})
under the name of ``canonical basis'' for the A,D,E cases.
The dual basis of the lower crystal basis of $\Uf$ is called
the {\em upper global basis} of $\Uf$.

\subsection{Symmetry}
Let $\theta$ be an automorphism of
$I$ such that $\theta^2=\id$ and 
$(\al_{\theta(i)},\al_{\theta(j)})=(\al_i,\al_j)$.
Hence it extends to an automorphism of the root lattice $Q$
by $\theta(\al_i)=\al_{\theta(i)}$,
and induces an automorphism of $\U$.

Let $\B_\theta(\g)$ be the $K$-algebra
generated by $E_i$, $F_i$, and
invertible elements $T_i$ ($i\in I$)
satisfying the following defining relations:
\eq&&\left\{
\parbox{30em}{
\bi
\item the $T_i$'s commute with each other,
\item
$T_{\theta(i)}=T_i$ for any $i$,
\item
$T_iE_jT_i^{-1}=q^{(\al_i+\al_{\theta(i)},\al_j)}E_j$ and
$T_iF_jT_i^{-1}=q^{(\al_i+\al_{\theta(i)},-\al_j)}F_j$
for $i,j\in I$,
\item
$E_iF_j=q^{-(\al_i,\al_j)}F_jE_i+
(\delta_{i,j}+\delta_{\theta(i),j}T_i)$
for $i,j\in I$,
\item
the $E_i$'s and the $F_i$'s satisfy the Serre relations.
\ei}
\right.
\eneq
Hence $\B_\theta(\g)\simeq\Uf\otimes K[T_i^{\pm1};i\in I]\otimes \Ue$.
We set $E_i^{(n)}=E_i^n/[n]_i!$
and $F_i^{(n)}=F_i^n/[n]_i!$.

Let $\la\in P_+\seteq\set{\la\in \Hom(Q,\Q)}%
{\mbox{$\lan \al_i^\vee,\la\ran\in\Z_{\ge0}$ for any $i\in I$}}$ 
be a dominant integral weight such that
$\theta(\la)=\la$.
\Prop\label{prop:Vtheta}
\bi
\item
There exists a $\B_\theta(\g)$-module $V_\theta(\la)$
generated by a vector $\vac_\la$ such that
\be[{\rm(a)}]
\item
$E_i\vac_\la=0$ for any $i\in I$,
\item
$T_i\vac_\la=q^{(\al_i,\la)}\vac_\la$ for any $i\in I$,
\item
$\set{u\in V_\theta(\la)}{\text{$E_iu=0$ for any $i\in I$}}
=K\vac_\la$.
\ee
Moreover such a $V_\theta(\la)$ is irreducible and 
unique up to an isomorphism.
\item
there exists a unique symmetric bilinear form $(\scbul,\scbul)$
on $V_\theta(\la)$ such that $(\vac_\la,\vac_\la)=1$ and
$(E_iu,v)=(u,F_iv)$ for any $i\in I$ and $u,v\in V_\theta(\la)$,
and it is non-degenerate.
\ei
\enprop
The pair $(\Bt,V_\theta(\la))$
is an analogue of
$(\B,U_q^-(\g))$.
Such a $V_\theta(\la)$ is constructed as follows.
Let $\Uf\vac'_\la$ and $\Uf\vac''_\la$ be a
copy of a free $\Uf$-module.
We give the structure of a $\Bt$-module on them
as follows: for any $i\in I$ and $a\in\Uf$
\eq&&
\ba{rcl}
T_i(a\vac'_\la)&=&q^{(\al_i,\la)}(\Ad(t_it_{\theta(i)})a)\vac'_\la,\\[2pt]
E_i(a\vac'_\la)&=&
\bigl(e'_ia+q^{(\alpha_i,\la)}\Ad(t_i)(e^*_{\theta(i)}a)\bigr)\vac'_\la,\\[2pt]
F_i(a\vac'_\la)&=&(f_ia)\vac'_\la
\ea
\eneq
and
\eq&&\ba{rcl}
T_i(a\vac''_\la)&=&q^{(\al_i,\la)}(\Ad(t_it_{\theta(i)})a)\vac''_\la,\\[2pt]
E_i(a\vac''_\la)&=&(e'_ia)\vac''_\la,\\[2pt]
F_i(a\vac''_\la)&=&
\bigl(f_ia+q^{(\alpha_i,\la)}(\Ad(t_i)a)f_{\theta(i)}\bigr)\vac''_\la.
\ea
\eneq
Then there exists a unique $\Bt$-linear morphism $\psi\cl
\Uf\vac'_\la\to\Uf\vac''_\la$
sending $\vac'_\la$ to $\vac''_\la$.
Its image $\psi(\Uf\vac'_\la)$ is $V_\theta(\la)$.

Hereafter we assume further that
\eq&&
\text{there is no $i\in I$ such that $\theta(i)=i$.}
\eneq
We conjecture that $V_\theta(\la)$ has a crystal basis.
This means the following.
We define the modified root operators similarly to \eqref{eq:modified}:
$$\tE_i(u)=\sum_{n\ge1}F_i^{(n-1)}u_n\quad
\text{and}\quad
\tF_i(u)=\sum_{n\ge0}F_i^{(n+1)}u_n
$$
when writing $u=\sum_{n\ge0}F_i^{(n)}u_n$ with $E_iu_n=0$.
Let $L_\theta(\la)$ be the $\A_0$-submodule of
$V_\theta(\la)$ generated by $\tF_{i_1}\cdots\tF_{i_\ell}\vac_\la$
($\ell\ge0$ and $i_1,\ldots,i_\ell\in I$\,),
and let $B_\theta(\la)$ be the subset 
$\set{\tF_{i_1}\cdots\tF_{i_\ell}\vac_\la\bmod \qs L_\theta(\la)}%
{\text{$\ell\ge0$, $i_1,\ldots, i_\ell\in I$}}$ of
$L_\theta(\la)/\qs L_\theta(\la)$.
\begin{conjecture}
\bi
\item
$\tF_iL_\theta(\la)\subset L_\theta(\la)$
and $\tE_iL_\theta(\la)\subset L_\theta(\la)$,
\item
$B_\theta(\la)$ is a basis of $L_\theta(\la)/\qs L_\theta(\la)$,
\item
$\tF_iB_\theta(\la)\subset B_\theta(\la)$,
and
$\tE_iB_\theta(\la)\subset B_\theta(\la)\sqcup\{0\}$.
\ei
\end{conjecture}
Moreover we conjecture that
$V_\theta(\la)$ has a global crystal basis.
Namely, let $-$ be the bar-operator of
$V_\theta(\la)$ given by $-\cl a\vac_\la\to \bar a\vac_\la$ ($a\in \Uf$)
(such an operator exists).
\begin{conjecture}
$(L_\theta(\la),L_\theta(\la)^-,\Uf_\A\vac_\la)$
is balanced.
\end{conjecture}
Assume that this conjecture is
true. Let $G^\low\colon L_\theta(\la)/\qs L_\theta(\la)\isoto 
E\seteq L_\theta(\la)\cap L_\theta(\la)^-\cap \Uf_\A\vac_\la$
be the inverse of $E\isoto  L_\theta(\la)/\qs  L_\theta(\la)$.
Then $\set{G^\low(b)}{b\in B_\theta(\la)}$ forms a basis of $V_\theta(\la)$.
We call this basis the {\em lower global basis} of $V_\theta(\la)$.
Let $\set{G^\upper(b)}{b\in B_\theta(\la)}$ be the dual basis 
to $\set{G^\low(b)}{b\in B_\theta(\la)}$ with respect to the 
inner product of
$V_\theta(\la)$.
We call it the {\em upper global basis} of $V_\theta(\la)$.

We can prove the conjectures
in the $\gl_\infty$-case:
{\scriptsize
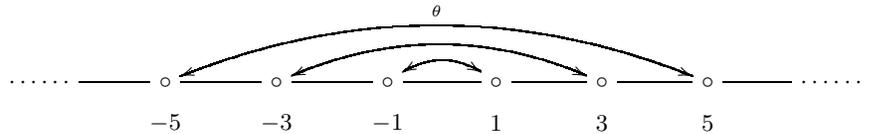
$$
\xymatrix@R=.8ex@C=6ex{
\cdots\cdots\ar@{-}[r]&\circ\ar@{-}[r]
\ar@/^1.8pc/@{<->}[rrrrr]^\theta&\circ\ar@{-}[r]\ar@/^1.2pc/@{<->}[rrr]&\circ\ar@{-}[r]
\ar@/^.7pc/@{<->}[r]&
\circ\ar@{-}[r]&\circ\ar@{-}[r]&\circ\ar@{-}[r]&\cdots\cdots\ .\\
&-5&-3&-1&\;1\;&\;3\;&\;5\;
}
$$}
\Th\label{th:main}
Let $I$ be the set $\Zodd$ of odd integers.
Define 
$$(\al_i,\al_j)=\begin{cases}
2&\text{if $i=j$,}\\
-1&\text{if $i=j\pm 2$,}\\
0&\text{otherwise,}
\end{cases}$$
and $\theta(i)=-i$.
Then, for $\la=0$,
$V_\theta(\la)$ has a crystal basis and a global basis.
\enth
Note that $\set{a\in \Uf}{a\vac_\la=0}=\sum_i\Uf(f_i-f_{\theta(i)})$
in this case.

The proof is by using a kind of PBW basis,
similarly to \cite{L}.
The details will appear elsewhere.

The following diagram is the part of the crystal graph of $B_\theta(\la)$
that concerns only the $1$-arrows and the $(-1)$-arrows.
$$
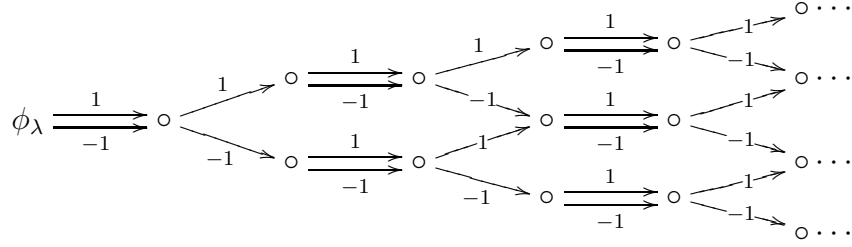
\xymatrix@R=.08em@C=3em{
&&&&&&\circ\cdots\\
&&&&\circ\ar@<.2pc>[r]^{1}\ar@<-.2pc>[r]_{-1}&
\circ\ar[ru]|{1}\ar[rd]|(.45){-1}\\
&&\circ\ar@<.2pc>[r]^{1}\ar@<-.2pc>[r]_{-1}&\circ
\ar[ru]^{1}\ar[rd]|{-1}&&&\circ\cdots&\\
\vac_\la\ar@<.2pc>[r]^{1}\ar@<-.2pc>[r]_{-1}&
\circ\ar[ru]^{1}\ar[rd]_{-1}&&&\circ\ar@<.2pc>[r]^{1}\ar@<-.2pc>[r]_{-1}
&\circ\ar[ru]|{1}\ar[rd]|(.45){-1}\\
&&\circ\ar@<.2pc>[r]^{1}\ar@<-.2pc>[r]_{-1}&\circ\ar[ru]|{1}\ar[rd]_{-1}
&&&\circ\cdots\\
&&&&\circ\ar@<.2pc>[r]^{1}\ar@<-.2pc>[r]_{-1}&\circ\ar[ru]|{1}
\ar[rd]|(.45){-1}
\\
&&&&&&\circ\cdots
}$$
Here is the part of the crystal graph of $B_\theta(\la)$
that concerns only the $n$-arrows and the $(-n)$-arrows
for an odd integer $n\ge 3$:
$$\xymatrix@C=8ex{
\vac_\la\ar@<.2pc>[r]^(.5){n}\ar@<-.2pc>[r]_(.45){-n}&
\circ\ar@<.2pc>[r]^(.45){n}\ar@<-.2pc>[r]_(.4){-n}&
\circ\ar@<.2pc>[r]^(.45){n}\ar@<-.2pc>[r]_(.4){-n}&
\circ\ar@<.2pc>[r]^(.4){n}\ar@<-.2pc>[r]_(.35){-n}&
\circ\cdots
}$$

\section{Affine Hecke algebra of type B}\label{sec:aff}
\subsection{Definition}
For $p_0,p_1\in\C^*$ and $n\in\Z_{\ge0}$, 
the affine Hecke algebra $\HB_n$ of type $B_n$
is the $\C$-algebra generated by
$T_i$ ($0\le i<n$) and invertible elements $X_i$ ($1\le i\le n$)
satisfying the defining relations:
\bi
\item
the $X_i$'s commute with each other,
\item
the $T_i$'s satisfy the braid relation:
$T_0T_1T_0T_1=T_1T_0T_1T_0$,
$T_iT_{i+1}T_i=T_{i+1}T_iT_{i+1}$ ($1\le i<n-1$),
$T_iT_j=T_jT_i$ ($\vert i-j\vert>1$),
\item
$(T_0-p_0)(T_0+p_0^{-1})=0$ and
$(T_i-p_1)(T_i+p_1^{-1})=0$ ($1\le i<n$).
\item
$T_0X_1^{-1}T_0=X_1$,
$T_iX_iT_i=X_{i+1}$ ($1\le i<n$),
and $T_iX_j=X_jT_i$ if $j\not=i,i+1$.
\ei
We assume that $p_0,p_1\in\C^*$ satisfy
\eq
p_0^2\not=1,p_1^2\not=1.
\eneq
Let us denote by $\Ft$ the Laurent polynomial ring
$\C[X_1^{\pm1},\ldots,X_n^{\pm1}]$,
and by $\Ftf$ its quotient field $\C(X_1,\ldots, X_n)$.
Then $\HB_n$ is isomorphic to the tensor product
of $\Ft$ and the subalgebra generated by the $T_i$'s
that is isomorphic to the Hecke algebra of type $B_n$.
We have
$$
T_ia=(s_ia)T_i+(p_i-p_i^{-1})\dfrac{a-s_ia}{1-X^{-\al_i^\vee}}
\quad\text{for $a\in\Ft$.}$$
Here $p_i=p_1$ ($1<i<n$),
and $X^{-\al_i^\vee}=X_1^{-2}$ ($i=0$)
and $X^{-\al_i^\vee}=X_{i}X_{i+1}^{-1}$ ($1\le i<n$).
The $s_i$'s are the Weyl group action on $\Ft$:
$(s_ia)(X_1,\ldots,X_n)=a(X_1^{-1},X_2,\ldots,X_n)$ for $i=0$
and
$(s_ia)(X_1,\ldots,X_n)=a(X_1,\ldots,X_{i+1},X_i,\ldots,X_n)$ for $1\le i<n$.

Note that $\HB_n=\C$ for $n=0$.

\subsection{Intertwiner}
The algebra $\HB_n$ acts faithfully on
$\HB_n/\sum_i\HB_n(T_i-p_i)\simeq\Ft$.
Set $\vphi_i=(1-X^{-\al_i^\vee})T_i-(p_i-p_i^{-1})\in\HB_n$
and $\tphi_i=(p_i^{-1}-p_iX^{-\al_i^\vee})^{-1}\vphi_i\in\Ftf\otimes_{\Ft}
\HB_n$.
Then the action of $\tphi_i$ on $\Ft$ coincides with $s_i$.
They are called intertwiners.

\subsection{Affine Hecke algebra of type A}\label{subsec:affA}
We will review the LLT conjecture,
reformulated and solved by S.{} Ariki,
on the affine Hecke algebras of type A.

The affine Hecke algebra $\HA_n$ of type $A_n$
is isomorphic to the subalgebra
of $\HB_n$ generated by $T_i$ ($1\le i<n$) and $X_i^{\pm1}$
($1\le i\le n$). For a finite-dimensional $\HA_n$-module $M$
let us decompose
\eq\label{eq:wtdec}
M=\soplus_{a\in(\C^*)^n}M_a
\eneq
where $M_a=\set{u\in M}{\text{$(X_i-a_i)^Nu=0$ for any $i$ 
and $N\gg0$}}$
for $a=(a_1,\ldots,a_n)\in(\C^*)^n$.
For a subset $I\subset\C^*$,
we say that $M$ is of type $I$ if all the eigenvalues of $X_i$
belong to $I$.
The group $\Z$ acts on $\C^*$ by $\Z\ni n\cl a\mapsto ap_1^{2n}$.
\Lemma
Let $I$ and $J$ be $\Z$-invariant subsets in $\C^*$ such that
$I\cap J=\emptyset$.
\bi
\item
If $M$ is an irreducible $\HA_m$-module of type $I$
and $N$ is an irreducible $\HA_n$-module of type $J$,
then $\Ind_{\HA_m\otimes\HA_n}^{\HA_{m+n}}(M\otimes N)$
is irreducible of type $I\cup J$.
\item
Conversely if $L$ is an irreducible $\HA_n$-module
of type $I\cup J$, then there exist $m$ $(0\le m\le n)$,
an irreducible $\HA_m$-module $M$ of type $I$
and an irreducible $\HA_{n-m}$-module $N$ of type $J$
such that
$L$ is isomorphic to $\Ind_{\HA_m\otimes\HA_{n-m}}^{\HA_{n}}(M\otimes N)$.
\ei
\enlemma
Hence in order to study the irreducible modules over the affine Hecke
algebras of type A,
it is enough to treat the irreducible modules of type $I$
for an orbit $I$ with respect to the $\Z$-action on $\C^*$.
Let $\KA_{I,n}$ be the Grothendieck group of
the abelian category of
finite-dimensional $\HA_n$-modules of type $I$.
We set $\KA_I=\soplus\nolimits_{n\ge0}\KA_{I,n}$.
Then $\KA_I$ has a structure of Hopf algebra 
where the product and the coproduct 
$$\mu\cl \KA_{I,m}\otimes \KA_{I,n}\to
\KA_{I,m+n},\quad
\Delta\cl\KA_{I,n}\to\soplus\limits_{i+j=n}\KA_{I,i}\otimes
\KA_{I,j}$$
are given by
$M\otimes N\mapsto\Ind_{\HA_m\otimes\HA_{n}}^{\HA_{m+n}}(M\otimes N)$
and
by $M\mapsto \Res^{\HA_n}_{\HA_i\otimes \HA_j}M$.
Let $\g_I$ be the Lie algebra associated to the Dynkin diagram
with $I$ as the set of vertices and with edges between $a$ and $ap_1^2$
($a\in I$). It means 
\eq
&&(\al_i,\al_j)=2\delta_{i,j}-\delta_{i,p_1^2j}-\delta_{p_1^2i,j}\quad
\text{for $i,j\in I$.}
\label{eq:Dy}
\eneq
Let $U_I$ be the unipotent group associated with the Lie subalgebra 
$\g_I^{-}$ of $\g_I$ generated
by the $f_i$'s.
Then we have
\Lemma Let $I$ be a $\Z$-invariant set.
Then $\C\otimes\KA_I$ is isomorphic to the algebra 
$\CO(U_I)$ of the regular functions on
$U_I$ as a Hopf algebra.
\enlemma
Here, for $a\in I$, $f_a$ corresponds to the one-dimensional $\HA_1$-module 
$\C_a$ on which $X_1$ acts by $a$.
Let $\{G^\upper(b)\}_{b\in B(\infty)}$
be the upper global basis of $\Uf$. Then
$\bigl(\soplus\C[q]G^{\upper}(b)\bigr)/\bigl((q-1)
\soplus \C[q]G^{\upper}(b)\bigr)$
is isomorphic to $\CO(U_I)$.
The following theorem is conjectured 
for Hecke algebras of type A
by Lascoux-Leclerc-Thibon (\cite{LLT})
and reformulated and proved by S. Ariki (\cite{Ariki})
for affine Hecke algebras of type A.
\Th
The elements of $\KA$ associated to irreducible $\HA$-modules correspond
to the upper global basis $G^\upper(b)$
by the isomorphism above.
\enth
Hence the irreducible modules are parametrized by
$B(\infty)$.
Grojnowski (\cite{Groj}) constructed the operators $\te_a$ and $\tf_a$
on $B(\infty)$ in terms of irreducible modules.
The operator $\te_a$ sends an irreducible $\HA_n$ module $M$
to a unique irreducible submodule of the $\HA_{n-1}$-module
$\set{u\in M}{(X_n-a)u=0}$.
The operator $\tf_a$ sends an irreducible $\HA_n$ module $M$
to a unique irreducible quotient of
the $\HA_{n+1}$-module
$\Ind_{\HA_n\otimes\HA_1}^{\HA_{n+1}}(M\otimes \C_a)$.

\subsection{Representations of affine Hecke algebras of type B}
For $n,m\ge 0$, set
$\Fc\seteq\C[X_1^{\pm1},\ldots,X_{n+m}^{\pm1},
D^{-1}]$
where
$$D\seteq\prod\limits_{1\le i\le n<j\le n+m}\hs{-2ex}
(X_i-p_1^2 X_j)(X_i-p_1^{-2}X_j)
(X_i-p_1^2 X_j^{-1})(X_i-p_1^{-2}X_j^{-1})(X_i-X_j).$$
Then we can embed
$\HB_m$ into $\HB_{n+m}\otimes_{\Ft[n+m]}\Fc$
by 
\eqn
&&T_0\mapsto \tphi_n\cdots\tphi_1T_0\tphi_1\cdots\tphi_n,\ 
T_i\mapsto T_{i+n}\quad(1\le i<m),\ 
X_i\mapsto X_{i+n}\quad(1\le i\le m).
\eneqn
Its image commute with $\HB_n\subset\HB_{n+m}$.
Hence $\HB_{n+m}\otimes_{\Ft[n+m]}\Fc$
is a right $\HB_n\otimes\HB_m$-module.
\Lemma
$\HA_{n+m}\otimes_{\HA_n\otimes\HA_m}\bl(\HB_n\otimes\HB_m
\br)\otimes_{\Ft[n+m]}{\Fc}\isoto
\HB_{n+m}\otimes_{\Ft[n+m]}\Fc$.
\enlemma
For a finite-dimensional $\HB_n$-module $M$,
we decompose $M$ as in \eqref{eq:wtdec}.
The semidirect product group
$\Z_2\times\Z=\{1,-1\}\times\Z$ acts on $\C^*$
by $(\epsilon,n)\cl a\mapsto a^\epsilon p_1^{2n}$.

Let $I$ and $J$ be $\Z_2\times\Z$-invariant subsets of
$\C^*$ such that $I\cap J=\emptyset$.
Then for an $\HB_n$-module $N$
of type $I$ and $\HB_m$-module $M$ of type $J$,
the action of $\Ft[n+m]$ on $N\otimes M$ extends to an action of $\Fc$.
We set 
$$N\diamond M\seteq
(\HB_{n+m}\otimes_{\Ft[n+m]}\Fc)
\mathop\otimes\limits
_{(\HB_n\otimes\HB_m)\otimes_{\Ft[n+m]}{\Fc}}(N\otimes M).$$
By the lemma above, $N\diamond M$ is isomorphic to
$\Ind_{\HA_n\otimes\HA_m}^{\HA_{n+m}}(N\otimes M)$
as an $\HA_{n+m}$-module.
\Lemma
\bi
\item
Let $N$ be an irreducible $\HB_n$-module
of type $I$ and
$M$ an irreducible $\HB_m$-module of type $J$.
Then $N\diamond M$ is an irreducible $\HB_{n+m}$-module
of type $I\cup J$.
\item Conversely if $L$
is an irreducible $\HB_{n}$-module
of type $I\cup J$, then there exist an integer $m$ $(0\le m\le n)$,
an irreducible $\HB_m$-module $N$
of type $I$ and
an irreducible $\HB_{n-m}$-module $M$ of type $J$
such that $L\simeq N\diamond M$.
\item
Assume that a $\Z_2\times\Z$-orbit $I$ decomposes into
$I=I_+\sqcup I_-$ where $I_\pm$ are $\Z$-orbits and $I_-=(I_+)^{-1}$.
Assume that $\pm1,\pm p\not\in I$.
Then for any irreducible $\HB_n$-module $L$ of type $I$,
there exists an irreducible $\HA_n$-module $M$ such that
$L\simeq\Ind_{\HA_n}^{\HB_n}M$.
\ei
\enlemma
Hence in order to study $\HB$-modules, it is enough to
study irreducible modules of type $I$ for a $\Z_2\times\Z$-orbit
$I$ in $\C^*$ such that $I$ is a $\Z$-orbit or $I$
contains one of $\pm1,\pm p$.

In this note, {\em we treat only the case when 
the $\Z_2\times\Z$-orbit $I$ does not contain $1$ nor $-1$.}

For a  $\Z_2\times\Z$-invariant subset $I$ of $\C^*$,
we define
$\KB_{I,n}$
and $\KB_I$ similarly to the case of $A$-type.
Then $\KB_I$ is a (right) Hopf $\KA_I$-bimodule
by the multiplication and the comultiplication
$$\mu\cl\KB_{I,n}\times \KA_{I,m}\to \KB_{I,n+m}\quad\text{and}\quad
\Delta\cl\KB_{I,n}\to\soplus_{i+j=n}\KB_{I,i}\otimes\KA_{I,j}$$
given by $L\otimes M\mapsto \Ind_{\HB_n\otimes\HA_m}^{\HB_{n+m}}(L\otimes M)$
and $L\mapsto\Res^{\HB_n}_{\HB_i\otimes\HA_j}L$.

Let $\theta$ be the automorphism of $I$ given by $a\mapsto a^{-1}$.
Then it induces an automorphism of $U_I$.
Let $U_I^\theta$ be the $\theta$-fixed point sets of $U_I$.
Then the action of $\CO(U_I)\simeq\C\otimes\KA_I$
on $\KB_I$ descends to an action of $\CO(U_I^\theta)$,
as it follows from the following lemma.
\Lemma For
an irreducible $\HB_n$-module $L$ and
an irreducible
$\HA_m$-module $M$,
we have $\mu(L\otimes M)=\mu(L\otimes M^\theta)$,
where $M^\theta$ is the $\HA_m$-module induced from $M$
by the automorphism of $\HA_m$
given by $X_i\mapsto X_{m+1-i}^{-1}$, $T_i\mapsto T_{m-i}$.
\enlemma

\medskip
Now we take the case
$$I=\set{p_1^n}{n\in\Zodd}.$$
Assume that any of $\pm1$ and $\pm p_0$ 
is not contained in $I$.
The set $I$ may be regarded as the set of vertices of a Dynkin diagram
by \eqref{eq:Dy}.
Let us define an automorphism $\theta$ of $I$ by $a\mapsto a^{-1}$.
Let $\g_I$ be the associated Lie algebra
($\g_I$ is isomorphic to $\gl_\infty$ if $p_1$ has an infinite order,
and 
isomorphic to $A^{(1)}_\ell$
if $p_1^2$ is a primitive $\ell$-th root of unity).
Let $V_\theta(\la)$ be as in Proposition~\ref{prop:Vtheta} with $\la=0$.
\begin{conjecture}
\be[(i)]
\item
$V_\theta(\la)$ has a crystal basis and a global basis.
\label{(i)}
\item
the elements of $\KB_I$ associated with irreducible representations
corresponds to
the upper global basis of $V_\theta(\la)$ at $q=1$.
\ee
\end{conjecture}
Note that \eqref{(i)} is nothing but
Theorem~\ref{th:main} when $p_1$ is not a root of unity.

\medskip
Let us take the case
$$I=\set{p_0p_1^{2n}}{n\in\Z}\cup\set{p_0^{-1}p_1^{2n}}{n\in\Z}.$$
Assume that there exists no integer $n$ such that $p_0^2=p_1^{4n}$.
It includes the case where $p_0=p_1$ and $p_1^{2n}\not=1$ for any $n\in\Zodd$.
Let $\theta$ be the automorphism of $I$ given by
$\theta\cl a\mapsto a^{-1}$. Then $\theta$ has no fixed points.
We regard $I$ as the set of vertices of a Dynkin diagram by \eqref{eq:Dy}.
Let $\g_I$ be the associated Lie algebra.
It is isomorphic to either $\gl_\infty\oplus\gl_\infty$,
$\gl_\infty$, $A^{(1)}_\ell\oplus A^{(1)}_\ell$ or $A^{(1)}_\ell$.
Set $\la=\Lambda_{p_0}+\Lambda_{p_0^{-1}}$
(i.e.{} $(\al_i,\la)=\delta_{i,p_0}+\delta_{i,p_0^{-1}}$).

\begin{conjecture}
\bi
\item
$V_\theta(\la)$ has a crystal basis and a global basis.
\item
the elements of $\KB_I$ associated with irreducible representations
corresponds to
the upper global basis of $V_\theta(\la)$ at $q=1$.
\ei
\end{conjecture}
In the both cases, we conjecture that,
for an irreducible $\HB_n$-module $M$ corresponding to an upper global basis
$G^\upper(b)$, $\dim M_a$
coincides with the value of 
$(\vac_\la, E_{a_1}\cdots E_{a_n}G^\upper(b))$ at $q=1$
for $a=(a_1,\ldots,a_n)\in I^n$.

Miemietz (\cite{Mie}) introduced 
the operators $\tilde e_i$ and $\tilde f_i$ 
on the set of isomorphic classes of
irreducible modules, similarly to the A type case,
and studied their properties.
We conjecture that they coincide with
the operators $\tE_i$ and $\tF_i$ on $B_\theta(\la)$.

{\em Acknowledgment\quad}
The authors thank Syu Kato
for helpful discussions. They thank also
Bernard Leclerc for his comments on this paper.

\end{document}